\documentclass[a4paper]{article}
\usepackage{a4wide,amsmath,amsthm,amssymb}
\usepackage[english]{babel}
\usepackage{graphicx}
\usepackage[numbers]{natbib}
\usepackage{ifpdf}
\usepackage{hyperref}

\hypersetup{
    bookmarksopen=false,
    bookmarksnumbered=true,
    pdftitle={On open problems in polling systems},
    pdfauthor={M.A.A. Boon, O.J. Boxma, E.M.M. Winands}
}
\ifpdf
  \hypersetup{colorlinks=true,linkcolor=black,urlcolor=black,citecolor=black,pdfpagemode=UseOutlines,plainpages=false,pdfpagelabels}
\else
  \hypersetup{colorlinks=true,linkcolor=black,urlcolor=black,citecolor=black,pdfpagemode=UseOutlines,plainpages=false,pdfpagelabels}
\fi

\newtheorem{theorem}{Theorem}[section]

\newtheorem{conjecture}[theorem]{Conjecture}

\newtheorem{problem}{Open problem}
\theoremstyle{definition}

\newtheorem{remark}[theorem]{Remark}

\newcommand{\E}{\mathbb{E}}
\newcommand{\cp}{c_{\!P_1}^2}
\providecommand{\href}[2]{#2}

\title{On open problems in polling systems\footnote{The research was done in the framework of the BSIK/BRICKS project, and of the European Network of Excellence Euro-NF.}}
\author{Marko Boon\footnote{\textsc{Eurandom} and Department of Mathematics and Computer Science, Eindhoven University of Technology, P.O. Box 513, 5600MB Eindhoven, The Netherlands}\\\href{mailto:marko@win.tue.nl}{marko@win.tue.nl} \and O.J. Boxma\footnotemark[2]\\\href{mailto:boxma@win.tue.nl}{boxma@win.tue.nl} \and E.M.M. Winands\footnote{Department of Mathematics, Section Stochastics, VU University, De Boelelaan 1081a, 1081HV Amsterdam, The Netherlands}\\\href{mailto:emm.winands@few.vu.nl}{emm.winands@few.vu.nl}}

\date{January, 2011}
\begin{document}
\maketitle

\begin{abstract}
In the present paper we address two open problems concerning polling systems, viz., queueing systems consisting of multiple queues attended by
a single server that visits the queues one at a time. The first open problem deals with a system consisting of two queues, one of which has gated service, while the other receives 1-limited service. The second open problem concerns polling systems with general (renewal) arrivals and deterministic switch-over times that become infinitely large. We discuss related, known results for both problems, and the difficulties encountered when trying to solve them. %Nevertheless, our feeling is that both problems might be solved, which makes them very challenging.

\end{abstract}

\section{Introduction}

A polling system is a queueing system consisting of multiple queues attended by
a single server that
visits the queues one at a time. Polling systems naturally arise in a large
number of application areas, like
\begin{itemize}
\item
maintenance: a patrolling repairman visits various sites;
\item
manufacturing: a machine successively produces items of various types;
\item
computer-communication systems: a central computer cyclically polls the
terminals on a common link
to inquire whether they have any data to transmit;
\item
road traffic: traffic lights determine which traffic streams may proceed.
\end{itemize}
In many of these applications, the server incurs a non-negligible switch-over time when switching between queues.

There is a huge body of literature on polling systems, in which the basic cyclic polling system and many enhancements have been studied. Extensive surveys on polling systems and their applications may be found in \cite{levysidi90,takagi1988qap,vishnevskiisemenova06}. In this note we present two challenging open problems motivated by two of the aforementioned application areas. By doing so, we want to stimulate research and  new collaborations in these directions.
%, which in turn give birth to challenging new queueing theoretic problems.
The first problem is motivated by a computer-communication system application and seems to lead to a boundary value problem with a rather complicated shift. The second problem requires an asymptotic analysis of the waiting-time distribution and stems from a manufacturing application.

\section{Model and notation}

In the present paper (and in almost the whole polling literature), the server visits the $N$ queues in cyclic order $Q_1,Q_2,\dots,Q_N,Q_1,\dots$, and the arrival processes of customers at the various queues are assumed to be independent Poisson processes, with rate $\lambda_i$ at $Q_i$, $i=1,\dots,N$. The service requirements at $Q_i$, denoted by $B_i$, are independent, identically distributed random variables with LST (Laplace-Stieltjes transform) $\beta_i(\cdot)$, $i=1,\dots,N$. Similarly, the switch-over times between $Q_i$ and $Q_{i+1}$, denoted by $S_i$, are independent, identically distributed with LST $\sigma_i(\cdot)$, $i=1,\dots,N$. The sum of the switch-over times is denoted by $S$. We furthermore assume that all arrival, service and switch-over processes are independent, and that the various parameters are such that the joint steady-state queue-length distributions at server visit epochs, server departure epochs and arbitrary epochs exist. We also introduce the notation $\rho_i = \lambda_i\E[B_i]$, and $\rho = \sum_{i=1}^N \rho_i$.

A key aspect of polling systems is the service discipline at each queue. The three most important service disciplines are {\em exhaustive} (E): a queues is served until it is empty; {\em gated} (G): the server only serves those customers which were present at the start of the visit; and {\em 1-limited} ($1$-L): the server serves only one customer -- if any is present. In a seminal paper, Resing \cite{Resing} has shown that the PGF (Probability Generating Function) of the joint steady-state queue-length distribution at epochs at which the server arrives at, say, $Q_1$ can be obtained explicitly for those polling systems in which the service discipline at each queue is of a branching-type, viz., the following holds for $i=1,\dots,N$:
\\
{\em If there are $k_i$ customers present at $Q_i$ at the start of a visit, then during the course of the visit,
each of these $k_i$ customers will effectively be replaced
in an i.i.d. manner by a random population having PGF
$h_i(z_1,\dots,z_N)$, which may be any $N$-dimensional PGF.}
\\
The joint queue-length process at visit epochs then becomes an $N$-class branching process
with immigration (the immigration corresponding to arrivals during switch-over times).
One may easily verify that exhaustive and gated are branching-type disciplines, whereas $1$-limited is not.
In the gated case, $h_i(z_1,\dots,z_N) = \beta_i(\sum_{j=1}^N \lambda_j(1-z_j))$,
and in the exhaustive case,
$h_i(z_1,\dots,z_N) = \pi_i(\sum_{j \neq i}^N \lambda_j(1-z_j))$,
where $\pi_i(\cdot)$ denotes the LST of the busy period distribution at $Q_i$, when viewed
as an $M/G/1$ queue in isolation.

A key element in the analysis of such branching-type polling systems is that the following relation
holds between the PGF $G_i(z_1,\dots,z_N)$ of the joint queue-length distribution
at the {\em end} of a visit to $Q_i$ and the PGF $F_i(z_1,\dots,z_N)$ of the joint queue-length
distribution at the {\em start} of that visit:
\begin{equation}
G_i(z_1,\dots,z_N) = F_i(z_1,\dots,z_{i-1}, h_i(z_1,\dots,z_N),z_{i+1},\dots,z_N).
\label{g}
\end{equation}
Moreover, it is easily seen that
\begin{equation}
F_{i+1}(z_1,\dots,z_N) = \sigma_i\big(\sum_{j=1}^N \lambda_j(1-z_j)\big) G_i(z_1,\dots,z_N).
\label{f}
\end{equation}
Successively applying each of these equations once for $i=1,\dots,N$,
one may now express $F_1 (z_1,\dots,z_N)$ into itself.
After iteration this yields an expression for $F_1(z_1,\dots,z_N)$ in the form of an infinite sum of products.

In the next two sections we formulate two open problems for polling systems
with a partial, respectively full, branching-type service discipline.

\section{\boldmath Open problem 1: gated plus $1$-limited}
In this section we restrict ourselves to the case of $N=2$ queues.
We are interested in determining $F_1(z_1,z_2)$ and $F_2(z_1,z_2)$.
After briefly discussing known results in Subsection~\ref{2.1}, we formulate in Subsection~\ref{2.2} an open problem
regarding the polling model with a gated and a 1-limited queue.

\subsection{Known results for two-queue polling systems}
\label{2.1}
The polling models with discipline E/E (Exhaustive/Exhaustive), G/G and E/G
fall in the class of multi-type branching, and are easily solved;
E/E was already solved by Tak\'acs \cite{Takacs} in 1968.
We refer to the survey \cite{takagi1988qap} and to \cite{Resing} for the other cases and for extensions
to a general number of queues $N$.
The $1$-L/$1$-L model was solved by using the theory of boundary value problems;
see \cite{CohenBoxma} for the case of zero switch-over times, and \cite{BoxmaGroenendijk}
for the case of non-zero switch-over times.
Now let us turn to E/1-L and G/1-L.
First observe that, with $Q_2$ having 1-L:
\begin{equation}
G_2(z_1,z_2) =  \frac{\beta_2(\lambda_1(1-z_1) + \lambda_2(1-z_2))}{z_2} [F_2(z_1,z_2) - F_2(z_1,0)]
+ F_2(z_1,0).
\label{1L}
\end{equation}
Combination of \eqref{g}, \eqref{f} and \eqref{1L} yields,
after having introduced $\beta_i(z_1,z_2) := \beta_i(\lambda_1(1-z_1)+\lambda_2(1-z_2))$
and
$\sigma_i(z_1,z_2) := \sigma_i(\lambda_1(1-z_1)+\lambda_2(1-z_2))$, $i=1,2$:
\begin{eqnarray}
F_{1}(z_1,z_2) &=& \frac{\beta_2(z_1,z_2) \sigma_2(z_1,z_2)}{z_2}
[\sigma_1(z_1,z_2) F_1(h_1(z_1,z_2),z_2) - \sigma_1(z_1,0) F_1(h_1(z_1,0),0)]
\nonumber
\\
&+& \sigma_2(z_1,z_2) \sigma_1(z_1,0) F_1(h_1(z_1,0),0).
\label{branch1-L}
\end{eqnarray}
The E/$1$-L model with zero switch-over times is simply a two-class nonpreemptive priority
model.
Ibe \cite{Ibe} considers the case of non-zero switch-over times, obtaining
the marginal queue-length distribution in $Q_1$ at polling instants of that queue.
It is less well-known that the {\em joint} queue-length distributions at polling instants of a queue
can also be found in a quite straightforward manner.
This is accomplished by substituting $h_1(z_1,z_2) = \pi_1(\lambda_2(1-z_2))$
into (\ref{branch1-L}), calling this function $g(z_2)$, and observing that $F_1(h_1(z_1,0),0)=F_1(g(0),0)$
is a constant, say $C$, not depending on $z_1$:
\begin{eqnarray}
F_{1}(z_1,z_2) &=& \frac{\beta_2(z_1,z_2) \sigma_2(z_1,z_2)}{z_2}
[\sigma_1(z_1,z_2) F_1(g(z_2),z_2) - C \sigma_1(z_1,0)]
\nonumber
\\
&+& C \sigma_2(z_1,z_2) \sigma_1(z_1,0) .
\label{branch1-LA}
\end{eqnarray}
The substitution $z_1 = g(z_2)$ finally solves the problem.
Details may be found in Section 6.3 of the PhD thesis of Groenendijk \cite{Groenendijk}.
E/1-L appears to be conceptually easier than E/E or any other known polling model,
not requiring a branching-type sum-of-infinite-products solution, and neither the solution of a boundary value problem.

\begin{remark}
Our sketch of the analysis of E/1-L reveals that one can extend that analysis
to the case in which $h_1(z_1,z_2)$, which above equals $\pi_1(\lambda_2(1-z_2))$, is some arbitrary
PGF $g(z_2)$.
For example, when the server is at $Q_1$, one could have Poisson arrivals with rate $\lambda_2^*$, or batch Poisson arrivals, at $Q_2$.
\end{remark}

\begin{remark}
%The $2$-dimensional queue-length process in the $2$-queue polling model with exhaustive service at $Q_1$ and $k$-limited service at $Q_2$
%has been analyzed in {\em Erik, PEIS artikel noemen. Ook Ozawa, Lee?}.
For general $k$, an exact evaluation for the queue-length distribution is known in two-queue exhaustive/$k$-limited systems with zero setup times (see Lee \cite{lee96}
and Ozawa \cite{ozawa90,ozawa97}) and with state-dependent switch-over times (see \cite{winands09}).
\end{remark}

\begin{remark}
As discussed in various polling studies, one can readily derive the waiting-time LST at $Q_i$
from $F_i(z_1,z_2)$. Furthermore, there exists a simple relation between the {\em mean} waiting times $\E[W_i]$ (at $Q_i$), $i=1,2,\dots,N$
in polling systems, a so-called pseudo-conservation law (cf. \cite{Boxma89}).
In the G/1-L case, this pseudo-conservation law reduces to:
\begin{equation}
\rho_1 \E[W_1] + \rho_2\left(1 - \frac{\lambda_2 \E[S]}{1-\rho}\right) \E[W_2] =
\rho \sum_{i=1}^2 \frac{\lambda_i \E[B_i^2]}{2(1-\rho)} + \rho \frac{\E[S^2]}{2\E[S]} +
\frac{\E[S]}{2(1-\rho)}[\rho^2 + \rho_1^2 + \rho_2^2],
\end{equation}
\end{remark}

\subsection{An unsolved two-queue polling model: G/1-L}
\label{2.2}

Throughout the polling literature, G and E seem to have comparable complexity.
Their role in the branching-type polling models, and in (pseudo-)conservation laws, is similar.
In view of this, and of the simplicity of E/1-L,
it is remarkable that G/1-L has remained unsolved for the past twenty years,
despite the fact that it is a quite relevant model (cf.\ Bisdikian \cite{Bisdikian},
who introduces
a variant of
G/1-L %{\em ik zie geen verschil met G/1-L, maar mischien wil Marko ook nog een keer kijken}
as a model for communication networks with bridge-stations and suggests an approximative approach).
Hence we state

\begin{problem}
Determine the joint queue-length PGF at polling instants in the two-queue G/1-L polling system.
\end{problem}

In the G/1-L case, (\ref{branch1-L}) holds with $h_1(z_1,z_2) = \beta_1(z_1,z_2)$.
Obvious attempts to obtain $F_1(z_1,z_2)$ from (\ref{branch1-L}) include
substituting $z_1=0$ (which yields a derivative), and substituting $z_1 = h_1(z_1,z_2)$
into the lefthand side of (\ref{branch1-L}). In the latter case, one obtains terms
$F_1(h_1(h_1(z_1,z_2),z_2),z_2)$ and $F_1(h_1(h_1(z_1,z_2),0),0)$
in the righthand side, and iteration does not seem to lead to a solution.

Our feeling is that (\ref{branch1-L}), in combination with obvious analyticity conditions
of $F_1(\cdot,\cdot)$ inside the product of unit circles, leads to a boundary value problem
(cf.\ \cite{CB}),
but one with a rather complicated shift introduced by the function $h_1(\cdot,\cdot)$.
Boundary value problems with a shift have been studied in the Riemann-Hilbert framework
(cf.\ \cite{Gakhov}, Section 17), and even in the setting of polling systems (cf. \cite{Lee}, which studies a two-queue
polling model with Bernoulli service at both queues),
but the present problem seems particularly challenging.

\section{Open problem 2: switch-over time asymptotics}

The next open problem considers the class of polling systems, with $N$ queues, that allow a multi-type branching process interpretation. %This class contains the classical exhaustive and gated policies as special cases.
We are interested in the behaviour of the polling system (under proper scaling conditions), when the \textit{deterministic} switch-over times tend to infinity.  This large switch-over time problem is relevant from a practical point of view, since systems with large switch-over times find a wide variety of applications in manufacturing environments (see \cite{winandsPhD}).
Firstly, Subsection \ref{known2} summarises known results for switch-over time asymptotics in polling systems with Poisson arrivals. Subsequently, we pose a conjecture for the behaviour of systems in which the arrival process at each of the queues is a general (renewal) process (see Subsection \ref{conjecture2}). Finally, %in Subsection \ref{open2}
the rigourous proof of this conjecture is stated as an open problem.

\subsection{Known result for Poisson arrival processes} \label{known2}

\textit{Under the assumption of Poisson arrival processes}, Winands \cite{winands2} presents an exact asymptotic analysis of the waiting-time distribution in branching-type polling systems with \textit{deterministic} switch-over times when the switch-over times tend to infinity. The results of \cite{winands2} generalise  those derived in \cite{olsen,mei2,winands} for the special case of exhaustive and gated service. Since the waiting time grows to infinity in the limiting case, \cite{winands2} focusses on the asymptotic scaled waiting time $W_i/S$ as $S \rightarrow \infty$ while keeping the ratios of the switch-over times constant.
We introduce $\Phi_i$ as the ``exhaustiveness'' of the service discipline in $Q_i$, defined as $\Phi_i = 1-\frac{\partial}{\partial z_i}h_i(z_1,\dots,z_N)\big|_{z_1=1,\dots,z_N=1}$. Its interpretation is, that each customer present at the start of a visit to $Q_i$ will be replaced by a number of type $i$ customers with mean $1 - \Phi_i$. If $Q_i$ receives exhaustive service, the exhaustiveness is 1; for gated service, it is $1-\rho_i$. In case of Poisson arrivals and deterministic switch-over times, the distribution of the asymptotic scaled waiting time is given by
\begin{equation}
\frac{W_i}{S} \xrightarrow{d} \frac{1-\rho_i}{1-\rho} U_i, \qquad (S \rightarrow \infty),
\end{equation}
where $U_i$ is uniformly distributed on $[\frac{1-\Phi_i}{\Phi_i},\frac{1}{\Phi_i}]$.

The closed-form expression of the scaled delay distribution has an intuitively appealing interpretation. That is, in the case of increasing deterministic switch-over times the polling system converges to a deterministic cyclic system with continuous deterministic service rates $1/\E[B_i]$ and continuous demand rates $\lambda_i$, $i=1,2,\ldots,N$, which reveals itself, for example, in the fact that the scaled number of customers at $Q_i$ at a polling instant of $Q_i$ becomes deterministic in the limit as shown in \cite{winands2}. This means that in the limit the customers arrive to the system and are served at constant rates with no statistical fluctuation whatsoever and that the scaled queue lengths can be seen as continuous quantities. Therefore, the uniform distribution emerging in the limiting theorems can be explained by the fact that it represents the position of the server in the cycle on arrival of a tagged customer.

\begin{remark}
In \cite{mei} it is shown that in heavy traffic (HT), i.e., if the load tends to one, the impact of higher moments of the switch-over times on the waiting-time distribution vanishes. Consequently, the scaled asymptotic waiting time depends on the marginal switch-over time distributions only through the first moment of the total switch-over time in a cycle. Building upon this observation, \cite{winands2} analyses the scaled asymptotic waiting time in branching-type polling systems with \textit{generally distributed switch-over times under heavy traffic} when the switch-over times tend to infinity.
The behaviour of the polling system then becomes deterministic, just like in polling systems with deterministic switch-over times, which are not necessarily operating in HT.
% \hfill $\Box$
\end{remark}

\subsection{Conjecture for general renewal arrival processes} \label{conjecture2}

Until now, we have assumed that the arrival processes are Poisson processes.
This assumption is used in \cite{winands2} to derive the asymptotics presented in the previous subsection,
building upon a result of \cite{borst} which derives a strong relation between the waiting-time distributions in models \textit{with} and \textit{without} switch-over times. This relation is established by relating the similarities in the offspring generating functions of the underlying branching processes and by expressing the differences between the underlying immigration functions.
These results for polling systems with \textit{finite} switch-over times are exploited and, subsequently, it is shown that significant simplifications result as the switch-over times tend to \textit{infinity}. Unfortunately, the techniques used throughout \cite{winands2} rely heavily on the Poisson assumption, and corresponding results for polling systems with general arrival processes are not known.
Taking a second look at the intuitive interpretation of the aforementioned results, one would expect that the Poisson assumption is not essential for this kind of behaviour.
That is, if the scaled number of arrivals during an intervisit period
becomes deterministic, then the length of the scaled visit period (generated by these arrivals) converges to a constant as
well. Since the intervisit period is the sum of individual visit periods and switch-over times,
based on strong law of large numbers arguments the scaled number of arrivals subsequently
indeed tends to become deterministic. This circular intuitive reasoning (ignoring the
interdependence between the visit periods) is independent of the precise
characteristics of the renewal arrival process chosen.
We now pose this statement as a conjecture (see, also, \cite{winands2}).

\begin{conjecture}\label{renewalconjecture}
A cyclic polling system with general (renewal) arrival processes converges to a deterministic cyclic system when the deterministic switch-over times tend to infinity.
\end{conjecture}

To numerically test this conjecture for general arrival processes, we have performed a couple of simulation experiments of exhaustive polling systems with general renewal arrivals. %based on the simulation code described in \cite{vuuren}.
In Table \ref{tableMarcel}, we show results for a symmetric polling system with $3$ queues, where the service times are exponential with mean $0.25$. Interarrival times have mean $1$ and the corresponding squared coefficient of variation (SCV), $c^2_{A_i}$, is varied between $0.25$, $0.5$, $1$ and $2$. In order to obtain a distribution for these interarrival times, we fit a phase-type distribution on the first two moments (cf., e.g., \cite{tijms}). For the cases where the SCV equals $1$, Poisson processes are used for the arrival processes in order to obtain exact results, and this case is included as benchmark.
\begin{table}[ht]
\begin{small}
\begin{center}
\begin{tabular}{|l|c|c|c|c|}
\hline
\multicolumn{5}{|c|}{Cyclic}\\
\hline \hline
\rule{0cm}{2.5ex} & $c^2_{A_i}=0.25$ & $c^2_{A_i}=0.5$ & $c^2_{A_i}=1$ & $c^2_{A_i}=2$  \\ \hline
$S_i=1$ & $0.121$ & $0.167$ & \em0.259 & $0.444$ \\
$S_i=10$ & $0.012$ & $0.017$ & \em0.026 & $0.044$ \\
%$S_i=50$ & $0.002$ & $0.003$ & $0.005$ & $0.009$ \\
$S_i=100$ & $0.001$ & $0.002$ & \em0.003 & $0.004$ \\
 \hline
\end{tabular}
\caption{Squared coefficient of variation of the scaled number of customers at $Q_1$ at a polling instant of $Q_1$.  Values in \emph{italic} are not obtained by simulation, but are computed analytically.}
\label{tableMarcel}
\end{center}
\end{small}
\end{table}

Table \ref{tableMarcel} shows $\cp$, the SCV of the scaled number of customers at $Q_1$ at a polling instant of $Q_1$ for varying values of the marginal switch-over times $S_i$ in a cycle. From Table \ref{tableMarcel}, we clearly see that the coefficient of variation approaches zero when the switch-over times tend to infinity. For polling systems with deterministic switch-over times and Poisson arrivals, it can actually be shown analytically that the SCV of the number of customers in a queue, at the beginning of a visit to this queue, is inversely proportional to the total switch-over time $S$.
% See approx4.nb
Table \ref{tableMarcel} seems to suggest that this also holds for other arrival processes.
It goes without saying that a highly variable arrival process has a negative impact on how ``fast'' the limiting behaviour is approached. Via Chebyshev's inequality (see, e.g., \cite{papoulis1}) we know that a random variable with zero variance follows a deterministic distribution and, therefore, this observation provides empirical evidence for the fact that  the scaled number of customers at $Q_1$ at a polling instant of $Q_1$ becomes deterministic. Therefore, it confirms the validity of our conjecture that the polling system converges to a deterministic cyclic system as the switch-over times increase to infinity.

We have run tests for asymmetric polling systems as well, %For the sake of brevity of the present paper, we do not list the full numerical results, but we only mention that all tests give the same limiting results for cyclic polling systems. Other visit disciplines, like \emph{switch to the longest queue} and \emph{switch to a random queue}, have been studied as well.
as shown in Table \ref{tableMarko}. The first three columns show the input parameters: the SCV of the interarrival time distributions,~$c_{A_i}^2$, the imbalance of the interarrival times, $I_{A}$, and the imbalance of the service times, $I_{B}$.
The imbalance is the ratio between the largest and the smallest mean interarrival/service time. The arrival rates and mean service times are chosen such that the differences $\lambda_{i}-\lambda_{i+1}$ and $\E[B_{i+1}]-\E[B_i]$ are kept constant for $i=1,\dots,N-1$.
Furthermore, we have chosen the normalisation constraint $\sum_{i=1}^N\lambda_i/N = 1$, implying that the actual arrival rates and mean service times (for fixed $\rho$) follow from the relation $\rho=\sum_{i=1}^N\lambda_i\E[B_i]$.
See \cite{boonapprox2011} for a more elaborate description and some examples of how the arrival rates and mean service times can be computed from this definition of imbalance.
The last two columns of Table \ref{tableMarko} contain the SCVs of the waiting times of customers in $Q_1$, $c_{W_1}^2$, and the SCVs of the numbers of customers at $Q_1$ at a polling instant, $\cp$, for a cyclic polling system with $\rho=0.75$ and deterministic switch-over times $S_i = 100$, for $i=1,2,3$. The SCVs of the waiting times approach the limiting value $\frac13$, which is the SCV of a uniform distribution, quite rapidly. Furthermore, $\cp$ becomes negligibly small, illustrating that the behaviour of the system becomes deterministic.
\begin{table}[ht]
\[
\begin{array}{|ccc||rr|}
\hline
\multicolumn{5}{|c|}{\textrm{Cyclic}} \\
\hline
\hline
\rule{0cm}{2.5ex}c_{A_i}^2 & I_{A} & I_{B} & c_{W_1}^2 & \cp \\
\hline
0.25 &1 &1&   0.335&  0.001    \\
0.25 &1 &3&   0.335&  0.002    \\
0.25 &3 &1&   0.334&  0.001    \\
0.25 &3 &3&   0.335&  0.001    \\
\hline
1    &1 &1&  \emph{0.335} &  \emph{0.003}  \\
1    &1 &3&  \emph{0.336} &  \emph{0.003}  \\
1    &3 &1&  \emph{0.335} &  \emph{0.002}  \\
1    &3 &3&  \emph{0.336} &  \emph{0.003}  \\
\hline
2    &1 &1&  0.336 &  0.004   \\
2    &1 &3&  0.337 &  0.005   \\
2    &3 &1&  0.336 &  0.004   \\
2    &3 &3&  0.337 &  0.004   \\
\hline
\end{array}
\]
\caption{Squared coefficient of variation of the waiting time and the number of customers at $Q_1$ at a polling instant of $Q_1$. Values in \emph{italic} are not obtained by simulation, but are computed analytically.}
\label{tableMarko}
\end{table}

%\subsection{Open problem} \label{open2}

Summarising, we state the second open problem for polling systems.
\setcounter{problem}{1}
\begin{problem}
%Provide a rigourous analysis of the limiting waiting time distribution (under proper scaling conditions) in branching-type polling systems with general (renewal) arrival processes and \textit{deterministic} switch-over times, when the switch-over times tend to infinity.
Provide a rigourous proof of Conjecture \ref{renewalconjecture}, which states that a cyclic polling system with general (renewal) arrival processes converges to a deterministic cyclic system when the deterministic switch-over times tend to infinity.
%\hfill $\Box$
\end{problem}

We wish to end the present paper with stating a related open problem given in \cite{olsen}. That is, \cite{olsen} shows via numerical testing that similar limit theorems as presented here carry over to systems with Poisson arrivals and \textit{dynamic} visit orders (i.e., there exists no pre-determined order in which the queues are served). This phenomenon is intuitively explained via heuristic strong law reasoning in \cite{olsen} and it is conjectured that the limit theorems hold \textit{so long as the switch-overs perform a regulating effect}. As an example, we show simulation results in Table \ref{tableMarko2} for the same symmetric polling systems as studied in Table \ref{tableMarcel}, but now the server switches to the longest queue at the end of a visit. We can see clearly, that the system becomes deterministic as well. A resulting open problem is, therefore, the classification of polling systems in terms of service, visit and scheduling disciplines, which exhibit the discussed behaviour.

\begin{table}[ht]
\begin{small}
\begin{center}
\begin{tabular}{|l|c|c|c|c|}
\hline
\multicolumn{5}{|c|}{Longest queue}\\
\hline \hline
\rule{0cm}{2.5ex} & $c^2_{A_i}=0.25$ & $c^2_{A_i}=0.5$ & $c^2_{A_i}=1$ & $c^2_{A_i}=2$  \\ \hline
$S_i=1$     & 0.125 & 0.170 & 0.254 & 0.434 \\
$S_i=10$    & 0.012 & 0.017 & 0.026 & 0.044 \\
$S_i=100$   & 0.001 & 0.002 & 0.003 & 0.004 \\
 \hline
\end{tabular}
\caption{Squared coefficient of variation of the scaled number of customers at $Q_1$ at a polling instant of $Q_1$.}
\label{tableMarko2}
\end{center}
\end{small}
\end{table}

\section{Conclusions}

There is a huge literature on polling systems, due to their great applicability in real-life situations. In this paper we have described two problems that have remained unsolved in the polling literature, despite their practical relevance, and despite the fact that seemingly minor adaptations of these problems can be solved explicitly. For the first problem, which is the exact analysis of a two-queue polling system with respectively gated and 1-limited service, we pinpoint the difficulties one runs into when applying standard techniques. The second problem is the analysis of a polling system with general renewal arrivals under the limiting situation where the (deterministic) switch-over times tend to infinity. For this problem we have posed a strong conjecture stating that the (known) results for Poisson arrivals carry over to the system with general renewal arrivals. By posing these open problems we hope to provide a motivation to search for alternative ways to study and hopefully even solve them.

\bibliographystyle{abbrvnat}
%\bibliography{open}

\end{document}